\documentclass[a4paper]{amsart}
\usepackage{graphicx}
\usepackage{amsmath, amsthm}
\usepackage{amssymb}
\usepackage{amsfonts}
\usepackage{harvard}

\newtheorem{thm}{Theorem}[section]
\newtheorem{cor}[thm]{Corollary}
\newtheorem{prop}[thm]{Proposition}
\newtheorem{lem}[thm]{Lemma}
\numberwithin{equation}{section}

\newcommand\be{\begin{equation}}
\newcommand{\conv}{\operatorname{conv}}
\newcommand\ee{\end{equation}}

\newcommand{\var}{\operatorname{var}}

\newcommand\CC{{\mathbb C}}
\newcommand\EE{{\mathbb E}}

\newcommand\PP{{\mathbb P}}
\newcommand\RR{{\mathbb R}}

\begin{document}

\title{Functionals of Exponential Brownian Motion and Divided Differences}

\author{B. J. C. Baxter and R. Brummelhuis}

\address{School of Economics, Mathematics and Statistics,
Birkbeck College, University of London,
Malet Street, London WC1E 7HX, 
England}

\email{b.baxter@bbk.ac.uk}

\begin{abstract}
We provide a surprising new application of classical approximation
theory to a fundamental asset-pricing model of mathematical finance.
Specifically, we 
calculate an analytic value for the correlation coefficient between
exponential Brownian motion and its time average, and we find the
use of divided differences greatly elucidates formulae, providing a
path to several new results. As applications,
we find that this correlation coefficient is always at least
$1/\sqrt{2}$ and, via the Hermite--Genocchi integral relation, 
demonstrate that all
moments of the time average are certain 
divided differences of the exponential function. We also prove that
these moments agree with the somewhat more complex formulae 
obtained by Oshanin and Yor.
\end{abstract}

\maketitle

\section{Introduction}
We begin with exponential, or geometric, Brownian motion, defined by
\be
S(t) = e^{(r-\frac{\sigma^2}{2})t +
  \sigma B(t)}, \qquad t \ge 0,
\label{a1}
\ee
where $B : [0,\infty) \to \RR$ is Brownian motion. In other words,
$B$ is a stochastic process, or {\em random function}, for which
$B(0) = 1$, its increments are independent, and, for $0 \le s < t$, the
increment $B(t) - B(s)$ is normally distributed with mean zero and
variance $t-s$. The basic properties of Brownian motion are
explained in Section 37 of \citeasnoun{billingsley}, while 
\citeasnoun{karatzas} is a comprehensive treatise.
At a more elementary level, \citeasnoun{norris} provides a lucid derivation
of the main properties of Brownian motion, whilst \citeasnoun{higham}
provides a more general introduction well-suited to the numerical analyst.

We shall study the {\em time average}
\be
A(T) := \frac{1}{T} \int_0^T S(t)\,dt
\label{a2}
\ee
using the calculus of divided differences, a fundamental
tool in approximation theory. We will, in particular, show that the
correlation coefficient between $A(T)$ and $S(T)$, the moments of
$A(T)$ and, more generally, joint moments of $S(T)$ and $A(T)$ can be
elegantly, and usefully, expressed in terms of divided differences of
the exponential function. Now the time average $A(T)$ has been
extensively studied in the literature of Asian options; see, for
instance, \citeasnoun{yor}, \citeasnoun{yorbook} and
\citeasnoun{oshanin}. However, we find that our use of divided
differences both simplify and elucidate formulae. In Section 2, we
derive the correlation coefficient for $S(T)$ and $A(T)$, finding that
it is always at least $1/\sqrt{2}$, thus explaining the relative high
correlation that is observational folklore in the financial
community.\footnote{I am particularly grateful to my friend and colleague Dr
Dirk Siegel, a fellow PhD student of Mike Powell, who made this observation
known to me when consulting for
IBM.} In Section 3, we demonstrate that the divided
differences occurring in the lower moments of $S(T)$ and $A(T)$
generalise to all moments, using the fact that the integral of an
exponential function over a simplex can be expressed, via the
Hermite--Gennocchi formula, as a certain divided difference of the
exponential function. In Section 4, we provide the divided difference
theory required by the paper; some of this can, of course, be found in
the dedicatee's excellent textbook \citeasnoun{powell}, but
derivations of the
Hermite--Gennochi and Leibniz formula are less easily
available. Therefore we have provided their brief derivations in the
hope that this will enhance the 
paper's use to both the mathematical
finance and the numerical analysis communities. Finally, in Section 5,
we use our divided difference approach to derive a recurrence relation
for the moments of $A(T)$.

We first observe the familiar result
\be
\EE S(T) = e^{(r-\sigma^2/2)T} \EE e^{\sigma T^{1/2} Z}
= e^{(r-\sigma^2/2)T} e^{\sigma^2 T/ 2} = e^{rT}.
\label{a4}
\ee
Here $Z$ denotes a generic $N(0,1)$ Gaussian random variable and we have
used the standard fact that
\be
\EE e^{\lambda Z}
= (2\pi)^{-1/2} \int_\RR e^{\lambda \tau} e^{-\tau^2/2}\,d\tau
= (2\pi)^{-1/2} \int_\RR e^{-\frac{1}{2}\{(\tau - \lambda)^2 -
  \lambda^2\}}\,d\tau
= e^{\lambda^2/2}.
\label{a5}
\ee
Similarly,
\begin{eqnarray}
\EE A(T)
&=& T^{-1} \int_0^T \EE S(t)\,dt\nonumber\\
&=& \frac{e^{rT} - 1}{rT}.
\label{a6}
\end{eqnarray}
The approximation theorist will immediately recognise the divided difference
\be
\EE A(T) = \exp[0, rT],
\label{a7}
\ee
but a sceptical reader might view this as mere coincidence; in fact,
it is but the tip of an iceberg. 
We remind the reader that
$f[a_0, a_1, \ldots, a_n]$ is the highest coefficient of the unique polynomial
of degree $n$ interpolating $f$ at distinct points 
$a_0, \ldots, a_n \in \RR$, which implies $f[a_0] = f(a_0)$ and
\[
f[a_0, a_1] = \frac{f(a_1) - f(a_0)}{a_1 - a_0}.
\]
Further, it is evident that a divided difference does not depend
on the order in which the points $a_0, a_1, \ldots, a_n$ are chosen. 
As mentioned above, Section \ref{divdiffsec}
collects further divided difference theory required by this paper.

\section{The correlation coefficient between the time average and the asset}

We shall compute the correlation coefficient between $S(T)$
and $A(T)$. Specifically, we calculate
\be
R
:= \frac{\EE \left(S(T) A(T)\right) - 
      \EE \left(S(T)\right) \EE\left(A(T)\right)}
   {\sqrt{\var S(T) \var A(T)}}.
\label{a7_1}
\ee

\noindent
We find an elegant  divided difference expression for $R$.

\begin{thm} 
The correlation coefficient \eqref{a7_1} is given by
\be
R \equiv R(rT, \sigma^2 T) = \frac{ \exp[rT, 2rT, (2r+\sigma^2)T] }
    { \sqrt{ 2\exp[2rT, (2r+\sigma^2)T] \exp[0, rT, 2rT, (2r+\sigma^2)T]} }.
\label{a7_2}
\ee
\label{thm_corrcoeff}
\end{thm}

Let us begin our derivation.

\begin{lem} If $0 \le a \le b$, then
\be
\EE S(a) S(b) = \exp\Bigl( a (r+\sigma^2) + br \Bigr).
\label{a8}
\ee
\label{lem_a1}
\end{lem}

\begin{proof}
We have
\begin{eqnarray}
\EE S(a) S(b)
&=& \EE S(a)^2 e^{(b-a)(r-\sigma^2/2) + \sigma\left(B(b)-B(a)\right)}
\nonumber\\
&=& \EE S(a)^2 \EE e^{(b-a)(r-\sigma^2/2) + \sigma\sqrt{b-a} Z}\nonumber\\
&=& e^{(2r+\sigma^2) a} e^{(b-a)r}\nonumber\\
&=& e^{a(r+\sigma^2))} e^{br},
\end{eqnarray}
where $Z \sim N(0,1)$ and we have used \eqref{a4}.
\end{proof}

\begin{prop}
We have
\be
\EE S(T) A(T) = \exp[rT, (2r + \sigma^2)T].
\label{a9}
\ee
\end{prop}

\begin{proof}
Applying Lemma \ref{lem_a1}, we obtain
\begin{eqnarray*}
\EE S(T) A(T)
&=& T^{-1} \int_0^T \EE S(t) S(T)\,dt\\
&=& T^{-1} \int_0^T e^{(r+\sigma^2)t} e^{r T}\,dt\\
&=& \exp[rT, (2r+\sigma^2)T].
\label{a10}
\end{eqnarray*}
\end{proof}

\begin{prop} 
\be
\EE (A(T)^2) = 2 \exp[0, rT, (2r + \sigma^2)T].
\ee
\end{prop}

\begin{proof}
We find
\begin{eqnarray}
\EE (A(T)^2)
&=& T^{-2} \int_0^T \Bigl( \int_0^T \EE S(t_1) S(t_2)\, dt_2 \Bigr)\,dt_1 \nonumber\\
&=& 2T^{-2} \int_0^T \Bigl( \int_0^{t_1} \EE S(t_1) S(t_2)\, dt_2 \Bigr)\,dt_1.
\end{eqnarray}
Thus
\begin{eqnarray}
\EE (A(T)^2)
&=& 2T^{-2} \int_0^T \Bigl( \int_0^{t_1} e^{r(t_1 + t_2)} e^{\sigma^2 t_2}\, dt_2 \Bigr)\,dt_1 \nonumber\\
&=& 2T^{-2} \int_0^T e^{rt_1}\Bigl(\frac{e^{(r+\sigma^2)t_1} - 1}{r+\sigma^2} \Bigr)\,dt_1 \nonumber\\
&=& \frac{2}{(r+\sigma^2)T} \left[ \exp[0,(2r+\sigma^2)T] - \exp[0,rT] \right]\nonumber\\
&=& 2 \exp[0, rT, (2r+\sigma^2)T],
\end{eqnarray}
using the divided difference recurrence relation \eqref{a7b} to obtain the
final line.
\end{proof}

\noindent
Any reader still doubtful of the simplification provided by divided difference notation
might consider the alternative expression provided in \citeasnoun{hull}:
\[
\EE \left(A(T)^2\right) =
\frac{2 e^{(2r+\sigma^2)T}}{(r+\sigma^2)(2r+\sigma^2)T^2} + 
\frac{2}{rT^2}\left( \frac{1}{2r+\sigma^2} -
\frac{e^{rT}}{r+\sigma^2}\right).
\]
There is a similar divided difference relation for $\EE \left(A(T)^m\right)$, 
described
in the next section, but we now complete our derivation of Theorem \ref{thm_corrcoeff}.

\begin{proof}[Proof of Theorem \ref{thm_corrcoeff}]
Applying (\ref{a5}, \ref{a6}, \ref{a9}) and \eqref{a7b}, we obtain
\begin{eqnarray}
\EE S(T) A(T) - \EE S(T) \EE A(T)
&=& \exp[rT,(2r+\sigma^2)T] - e^{rT}(e^{rT} - 1)/(rT)\nonumber\\
&=& \exp[rT,(2r+\sigma^2)T] - \exp[rT,2rT]\nonumber\\
&=& \sigma^2 T \exp[rT, 2rT, (2r+\sigma^2)T].
\label{a15}
\end{eqnarray}

Further, 
\be
\var S(T) 
= \EE(S(T)^2) - (\EE S(T))^2
= e^{(2r+\sigma^2)T} - e^{2rT}
= \sigma^2 T \exp[2rT, (2r+\sigma^2)T],
\label{a16}
\ee
and, by (\ref{a6}, \ref{a10}), 
\begin{eqnarray}
\var A(T)
&=& 2\exp[0, rT, (2r+\sigma^2)T] - \left(\frac{e^{rT} - 1}{rT}\right)^2\nonumber\\
&=& 2\exp[0,rT,(2r+\sigma^2)T] - 2\exp[0,rT,2rT]\nonumber\\
&=& 2\sigma^2 T \exp[0,rT,2rT,(2r+\sigma^2)T],
\label{a17}
\end{eqnarray}
using the divided difference recurrence \eqref{a7b} once more.
Hence
\be
R = \frac{ \exp[rT, 2rT, (2r+\sigma^2)T] }
    { \sqrt{ 2\exp[2rT, (2r+\sigma^2)T] \exp[0, rT, 2rT, (2r+\sigma^2)T]} }.
\label{a18}
\ee
\end{proof}

\begin{figure}[tb]
\vspace{6pt}
\includegraphics[scale=0.8]{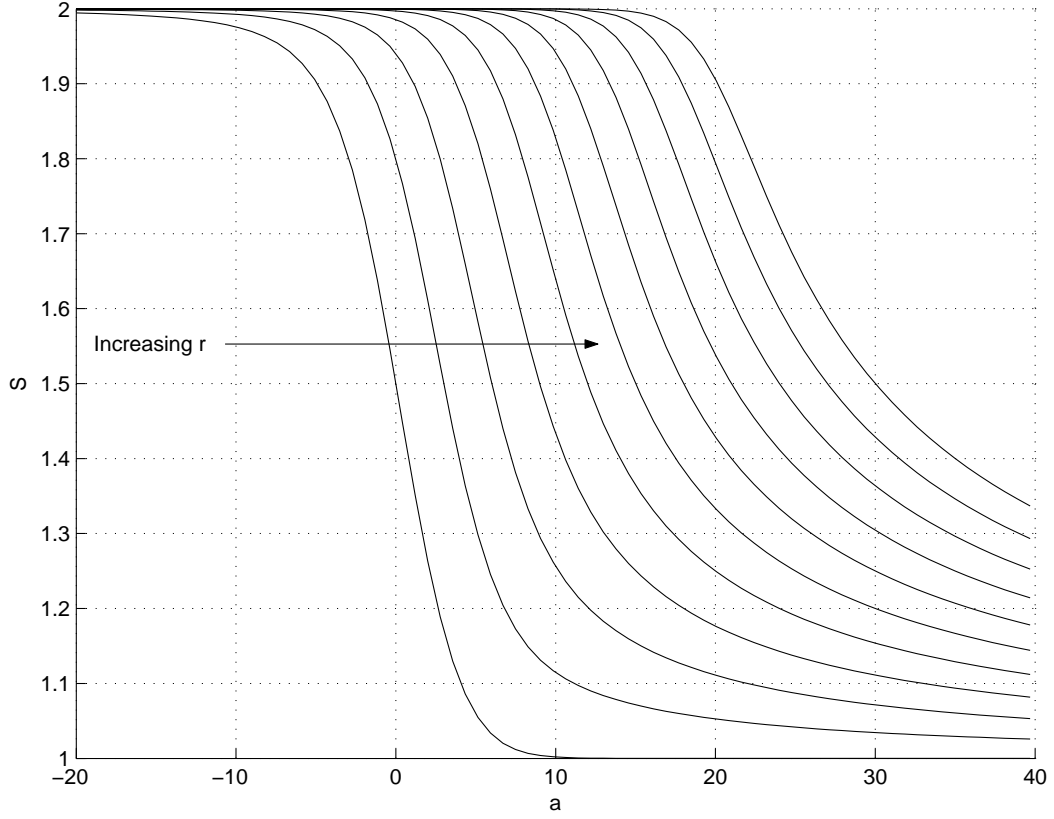}
\parbox{28.5pc}{\caption{$S(r,a)$ for $-20 \le a \le 40$ and $0.1 \le r
  \le 10$}\label{configs}}
\end{figure}

It is remarkable that the divided differences appearing in (\ref{a18})
are coefficients of the cubic polynomial interpolating the exponential
function at $0, rT, 2rT, (2r+\sigma^2)T$. We make three further observations:

\begin{enumerate}
\item

Armed with
an analytic expression for the correlation coefficient, we can apply
the exchange option valuation formula of \citeasnoun{margrabe}
to derive the values of certain Asian
options, if we are willing to accept that the time--average is suitably
approximated by exponential Brownian motion. We are investigating
the numerics of this rather simple approximation at present
and preliminary results are surprisingly promising.

\item
The correlation coefficient $R(rT, \sigma^2 T)$ is typically
close to unity: typical values of $r$, $\sigma$ and $T$ produce
values of $R$ in the $0.8-0.9$ range. In fact, we are able to prove
that the correlation coefficient satisfies $R(rT, \sigma^2 T) \ge
1/\sqrt{2}$, for all $r \ge 0$, $\sigma \ge 0$ and $T > 0$,
a surprisingly high lower bound for the correlation
coefficient. The details of this derivation are too complicated to
include here, and we refer the reader to \citeasnoun{baxter_fretwell} for
further details. However, the numerical findings are summarised in
Figure~\ref{configs}, which displays values of the closely related
quantity
\[
S \equiv S(r,a) = \frac{\left(\exp[a,2r,r]\right)^2}{\exp[a,2r]
  \exp[a,2r,r,0]},
\]
for $-20 \le a \le 40$ and $0.1 \le r \le 10$. It is easily checked
that $R(rT, \sigma^2 T) = \sqrt{S(rT, (2r+\sigma^2)T)/2}$, so that the
lower bound $R \ge 1/\sqrt{2}$ becomes $S \ge 1$. It is plausible that
$S(r,a)$ should be a decreasing function of $a$, for fixed $r$,
because correlation should be a decreasing function of
volatility. Further, it
is not difficult to establish the limiting values $\lim_{a \to -\infty}
S(a,r) = 2$ and $\lim_{r\to\infty} S(a,r) = 1$. However, further
analysis is not straightforward, and the analysis of
\citeasnoun{baxter_fretwell} makes great use of properties of
divided differences.

\item
It is natural to ask whether these divided difference expressions are
particular to exponential Brownian motion. In fact, similar
expressions occur when exponential Brownian motion is replaced by
certain L\'evy-stable variants: see \citeasnoun{BCF} for further details.

\end{enumerate}

\section{Computing higher moments of $A(T)$}

We now demonstrate that the neat divided difference formulae obtained for
the first and second moments of $A(T)$ are {\em not} coincidences,
but part of a greater pattern from which arise new formulae generalising
the moment calculations of \citeasnoun{oshanin}, \citeasnoun{yor} 
and \citeasnoun{yorbook}.

We begin with the iterated integral
\be
\EE A(T)^m
=
T^{-m} \int_0^T d\tau_m \int_0^T d\tau_{m-1} \cdots \int_0^T d\tau_1
\ \EE S(\tau_1) \cdots S(\tau_m).
\label{b3}
\ee
Now, given any point $(\tau_1, \ldots, \tau_m) \in [0,T]^m$,
let us sort its components into increasing order,
obtaining $(t_1, \ldots, t_n)$ (say). Then
\[
\EE S(\tau_1) \cdots S(\tau_m)
=
\EE S(t_1) \cdots S(t_m)
\]
and
\be
\EE A(T)^m
=
m! T^{-m}\int_0^T dt_m \int_0^{t_m} dt_{m-1} \cdots \int_0^{t_2} dt_1
\EE S(t_1) \cdots S(t_m).
\label{b4}
\ee
Our first task is to calculate the integrand, which we complete
after a simple lemma.

\begin{lem}
For any positive integer $k$, we have
\be
\EE \left[S(t)^k \right] = \exp \Bigl( krt + \frac{\sigma^2 t}{2} k(k-1)\Bigr).
\ee
\label{lem_b2}
\end{lem}

\begin{proof}
This is almost immediate from \eqref{a5}:
\[
\EE S(t)^k
= \EE e^{k(r-\sigma^2)t + \sigma k \sqrt{t} Z}
= e^{k(r-\sigma^2/2)t + \sigma^2 k^2 t/2}
= e^{krt + \sigma^2 t k(k-1)/2},
\]
where $Z \sim N(0,1)$.
\end{proof}

\begin{prop}
If $0 \le t_1 \le t_2 \le \cdots \le t_m$, then
\be
\EE S(t_1) S(t_2) \cdots S(t_m)
=
\exp \Bigl( \sum_{k=1}^m \left(r + (m-k)\sigma^2\right) t_k \Bigr).
\label{b5}
\ee
\label{prop_b3}
\end{prop}

\begin{proof} Lemma \ref{lem_a1} comprises the case $m=2$. We complete
the proof by induction on the number of terms $m$, 
first observing that, by a standard property of
geometric Brownian motion,
\be
\EE S(t_1) S(t_2) \cdots S(t_m)
= 
\EE S(t_1)^m \EE S(t_2 - t_1) \cdots S(t_m - t_1).
\label{b6}
\ee
Applying Lemma \ref{lem_b2} and our induction hypothesis, we obtain
\be
\EE S(t_1) S(t_2) \cdots S(t_m)
=
\exp\Bigl( mrt_1 + \sigma^2 t_1 m(m-1)/2
+ \sum_{\ell=2}^m 
\left( r + (m-\ell)\sigma^2\right) \left(t_\ell - t_1\right)\Bigr).
\label{b7}
\ee
The $t_1$ coefficient in the exponent is given by
\[
mr - (m-1)r + \sigma^2 t_1 \left( \frac{1}{2} m(m-1) - 
\sum_{\ell=1}^{m-2}\ell \right)
= r + \sigma^2 t_1 (m-1),
\]
using the elementary fact that $m(m-1)/2 = 1 + 2 + \cdots + m-1$. 
The coefficients of $t_2, \ldots, t_m$ are as already stated in \eqref{b5}.
\end{proof}

Thus the desired integral \eqref{b4} becomes
\begin{eqnarray}
\EE A(T)^m
&=& m! T^{-m} 
\int_0^T dt_m \int_0^{t_m} dt_{m-1} \cdots \int_0^{t_2} dt_1
\EE S(t_1) \cdots S(t_m)\nonumber\\
&=& 
m! \int_0^1 dt_m \int_0^{t_{m-1}} \cdots \int_0^{t_2} dt_1
\exp(\alpha_1 t_1 + \cdots \alpha_m t_m),
\label{b8}
\end{eqnarray}
where
\be
\alpha_k = \left(r + (m-k)\sigma^2 \right) T, \qquad k = 1, \ldots,
m.
\label{b9}
\ee

The integral displayed in \eqref{b8} can now be identified as a divided
difference using a variant form of the Hermite--Genocchi 
integral relation.

\begin{thm}
Let 
\be
b_k := k r + \sigma^2 k(k-1)/2, \qquad k = 0, 1, \ldots.
\label{b1}
\ee
Then
\be
\EE (A(T))^m
=
m! \exp[b_0 T, b_1 T, \ldots, b_m T], \qquad m \ge 0.
\label{b2}
\ee
\label{thm_b1}
\end{thm}

\begin{proof}
Apply Corollary \ref{dcor4} to \eqref{b8} and \eqref{b9}, using
$\sum_{k=1}^j k = j(j+1)/2$.
\end{proof}

The statement of Theorem \ref{thm_b1} simplifies when
$r = \sigma^2$, for then the drift term in \eqref{a1} vanishes, that is,
we consider $S(t) = \exp(\sigma \sqrt{t} B(t))$ alone; this
is the special case studied by \citeasnoun{oshanin} and
\citeasnoun{yor}, for the formulae grow {\em much} more complicated
without the use of divided differences. 
Therefore we now demonstrate that our expression
agrees with theirs.

\begin{thm}
If we set $r = \sigma^2/2$ in Theorem \ref{thm_b1}, then we obtain
\begin{align}
\EE \left[ A(T)^m \right]
&= m! \exp[0, rT, 2^2 rT, 3^2 rT, \ldots, m^2 rT]\nonumber\\
&= m! H_{\sqrt{rT}}[-m, \ldots, -1, 0, 1, \ldots, m],
\end{align}
where $H_{c}(x) := \exp(c^2 x^2)$, $x \in \RR$, for any positive $c$.
\label{thm_b2}
\end{thm}

\begin{proof}
We simply set $r = \sigma^2$ in Theorem \ref{thm_b1} and apply \eqref{dda0}.
\end{proof}

We can now apply Corollary \ref{dcor7} to derive the formula given in equation (14)
of \citeasnoun{oshanin}.

\begin{thm}
If we set $r= \sigma^2/2$, then
\be
\EE \left[A(T)^m \right]
=
\left( \frac{\Gamma(m)}{\Gamma(2m)}\right)
r^{-m} \left( -\frac{1}{2} (-1)^m {2m \choose m}
+ \sum_{\ell=0}^m {2m \choose \ell} (-1)^{\ell} e^{rT(m - \ell)^2} \right).
\label{osh1}
\ee
\end{thm}

\begin{proof}
Applying Corollary \ref{dcor7} to Theorem \ref{thm_b2}, we obtain
\begin{align}
\EE \left[A(T)^m\right]
&= \left( \frac{m!}{(2m)!} (rT)^{-m} \sum_{k=0}^{2m} {2m \choose k} 
(-1)^k e^{rT(k-m)^2}\right) \nonumber\\
&= 
\left( \frac{\Gamma(m)}{\Gamma(2m)}\right)
(rT)^{-m} \left( -\frac{1}{2} (-1)^m {2m \choose m}
+ \sum_{\ell=0}^m {2m \choose \ell} (-1)^{\ell} e^{rT(m - \ell)^2} \right),
\end{align}
after some straightforward algebraic manipulation.
\end{proof}

If we now replace $rT$ by $\alpha$ and $m$ by $j$ in \eqref{osh1}, then
we obtain equation (14) of \citeasnoun{oshanin}.

\section{Divided difference theory}
\label{divdiffsec}
Most of the properties of divided differences required here 
can be found in Chapter 5 of \citeasnoun{powell}. 
However,
proofs of the Hermite--Genocchi integral relation are
less easily available in the Anglophone mathematical literature, as is our
particular variant of it, although the specialist can find
much useful material in the treatise of \citeasnoun{lorentz}.
We have therefore provided a derivation for the convenience of the reader.
The Hermite--Genocchi formula and its consequences are
still very much topics of current research; see, for example, 
\citeasnoun{waldron}.
Furthermore, the result is better served in other European languages; see,
for instance, \citeasnoun{gelfond} for a French translation of a 
Russian classic,
or indeed the original \citeasnoun{hermite}.

We recall the {\em divided difference recurrence relation}.

\begin{thm}
\be
f[a_0, a_1, \ldots, a_n]
=
\frac{f[a_1, \ldots, a_n] - f[a_0, \ldots, a_{n-1}]}{a_n - a_0},
\label{a7b}
\ee
for any distinct complex numbers $a_0, \ldots, a_n$.
\end{thm}

\begin{proof}
See, for instance, \citeasnoun{powell}, Theorem 5.3.
\end{proof}

If $f$ is sufficiently differentiable, then we can, of course, define 
divided differences for
coincident points. Further, the elementary relation
\be
f[a_0, a_1] 
= \frac{f(a_1) - f(a_0)}{a_1 - a_0}
= \int_0^1 f^\prime((1-t)a_0 + ta_1)\,dt, \qquad\hbox{ when } a_0, a_1 \in \RR,
\label{dd1}
\ee
can be generalised to obtain the {\em Hermite--Genocchi formula}.

\begin{thm}[Hermite--Genocchi]
Let $f \in C^{(n)}(\RR)$ and let $a_0, a_1, \ldots, a_n$ be (not
necessarily distinct) real numbers
Then, for $n \ge 1$,
\begin{eqnarray}
\lefteqn{f[a_0, a_1, \ldots, a_n] }\nonumber \\
&=&
\int_{S_n} f^{(n)}(t_0 a_0 + t_1 a_1 + \cdots + t_n a_n)\,dt_1 \cdots
dt_n,\nonumber\\
&=& \int_0^1 dt_1 \int_0^{1-t_1} dt_2 \cdots
\int_0^{1 - \sum_{k=1}^{n-1}t_k} dt_n 
f^{(n)}(t_0 a_0 + t_1 a_1 + \cdots + t_n a_n)
\label{hg1}
\end{eqnarray}
where the domain of integration is the simplex
\be
S_n = \left\{ t = (t_1, t_2, \ldots, t_n) \in \RR_+^n \colon 
\sum_{k=1}^n t_k \le 1\right\}
\label{hg2}
\ee
and
\[
t_0 = 1 - \sum_{k=1}^n t_k.
\]
\end{thm}

\begin{proof}
We shall prove (\ref{hg1}) by induction on $n$, observing that
\[
\int_{S_1} f^\prime(t_0 a_0 + t_1 a_1)\,dt_1
=
\int_0^1 f^\prime(a_0 + t_1 (a_1 - a_0))\,dt_1
= \frac{f(a_1) - f(a_0)}{a_1 - a_0} = f[a_0, a_1].
\]
To extend the formula to higher-order divided differences, we 
note that
\be
f[a_0, a_1, \ldots, a_n, a_{n+1}]
=
\frac{f[a_1, a_2, \ldots, a_{n+1}] - f[a_0, a_1, \ldots, a_n]}{a_{n+1}
- a_0}
= g[a_0, a_{n+1}],
\label{hg3}
\ee
where
\be
g(x) = f[a_1, \ldots, a_n, x], \qquad x \in \RR.
\label{hg4}
\ee
Now
\[
g(x) = \int_{S_n} f^{(n)}(xt_0 + a_1 t_1 + \cdots + a_n t_n) \,dt_1
\cdots dt_n
\]
so that
\[
g^\prime(x) = \int_{S_n} t_0 f^{(n+1)}(xt_0 + a_1 t_1 + \cdots + a_n
t_n)\,dt_1 \cdots dt_n.
\]
Therefore
\begin{eqnarray*}
\lefteqn{f[a_0, a_1, \ldots, a_n, a_{n+1}]}\\
&=& \int_0^1 d\tau\ g^\prime( (1-\tau) a_0 + \tau a_{n+1} )\\
&=& \int_0^1 d\tau \int_{S_n} dt_1 \cdots dt_n\ t_0 f^{(n+1)}([(1-\tau)a_0
+ \tau a_{n+1}]t_0 + a_1 t_1 + \cdots a_n t_n)\\
&=& \int_{S_n} dt_1 \cdots dt_n \int_0^1 d\tau\ t_0
f^{(n+1)}([(1-\tau)t_0 a_0 + \sum_{\ell=1}^n a_\ell t_\ell + \tau t_0
a_{n+1})\\
&=& \int_0^1 dt_1 \int_0^{1-t_1} dt_2 \cdots
\int_0^{1-\sum_{k=1}^{n+1}t_k}
dt_{n+1} f^{(n+1)}\left(T_0 a_0 + \sum_{k=1}^{n+1} t_k a_k\right)\\
&=& \int_{S_{n+1}} f^{(n+1)}\left(T_0 a_0 + t_1 a_1 + \cdots + t_{n+1}
a_{n+1}\right) \,dt_1 \cdots dt_{n+1},\\
\end{eqnarray*}
where we have used the substitution $t_{n+1} = t_0 \tau$ and the
notation $T_0 = 1 - \sum_{k=1}^{n+1} t_k$.
\end{proof}

We shall need a variant form of the Hermite--Genocchi integral
relation for which the following notation is useful. Given
any real $n \times n$ nonsingular matrix $V$, with columns $v_1,
\ldots, v_n$, we let $K(V)$ denote 
the closed convex hull of 
$0, v_1, \ldots, v_n$, i.e.
\[
K(V) := \conv\{0, v_1, \ldots, v_n\}.
\]
In this notation, the Hermite--Gnocchi integral relation states that
\be
f[a_0, a_1, \ldots, a_n]
=
\int_{K(I_n)} f^{(n)}\left(a_0 + \left(a - a_0 e\right)^T
y\right)\,dy,
\ee
where 
\[
a = \begin{pmatrix} a_1 \\ \vdots \\ a_n \end{pmatrix}
,\qquad
e = \begin{pmatrix} 1 \\ \vdots \\ 1 \end{pmatrix},
\]
and $I_n$ denotes the $n \times n$ identity matrix.
Integrating the $n$th derivative over the simplex $K(V)$ yields
a useful variant form of Hermite--Genocchi.

\begin{thm} Let $V \in \RR^{n \times n}$ be any nonsingular matrix.
Then
\be
\frac{1}{\left|\det V\right|}
\int_{K(V)} f^{(n)}\left(a^T y\right)\,dy
=
f[0, (V^T a)_1, \ldots, (V^T a)_n],
\ee
where $(V^T a)_k$ denotes the $k$th component of the vector $V^T a$.
\label{dthm2}\end{thm}

\begin{proof}
Substituting $y = Vz$, Hermite--Genocchi implies the relation
\[
\int_{K(V)} f^{(n)}\left((V^T a)^T z\right)\,dz
= 
f[0, (V^T a)_1, \ldots, (V^T a)_n].
\]
\end{proof}

\begin{cor}
For any function $f \in C^{(n)}(\RR)$, we have
\begin{align}
&\int_0^1 dx_n \int_0^{x_n} dx_{n-1} \cdots \int_0^{x_2} dx_1
f^{(n)} \left(\sum_{k=1}^n a_k x_k\right)\nonumber\\
&\qquad =
f[0, a_n, a_n + a_{n-1}, \ldots, a_n + a_{n-1} + \cdots + a_1].
\end{align}
\label{dcor3}
\end{cor}

\begin{proof}
Set
\[
V = 
\begin{pmatrix}
1 & & & \\
1 & 1 & & \\
\vdots & & \ddots & \\
1 & 1 & \cdots & 1
\end{pmatrix}
\]
in Theorem \ref{dthm2}.
\end{proof}

The exponential function is a
particularly important case for us, 
in which case the Hermite--Genocchi formula becomes
\be
\exp[a_0, a_1, \ldots, a_n] = 
\int_{S_n} e^{t_0 a_0 + t_1 a_1 + \cdots + t_n a_n}\,dt_1 \cdots dt_n
\label{a7c}
\ee
and Corollary \ref{dcor3} takes the following form.

\begin{cor}
We have
\begin{align}
&\int_0^1 dx_n \int_0^{x_n} dx_{n-1} \cdots \int_0^{x_2} dx_1
\exp \left(\sum_{k=1}^n a_k x_k\right)\nonumber\\
&\qquad =
\exp[0, a_n, a_n + a_{n-1}, \ldots, a_n + a_{n-1} + \cdots + a_1].
\label{a7ca}
\end{align}
\label{dcor4}
\end{cor}

\begin{proof} 
Let $f$ be the exponential function in Corollary \ref{dcor3}.
\end{proof}

\noindent
Further, we note that, for the exponential function, Theorem \ref{dthm2}
becomes the interesting formula
\be
\frac{1}{\left|\det V\right|}
\int_{K(V)} e^{a^T y}\,dy
=
\exp[0, (V^T a)_1, \ldots, (V^T a)_n].
\ee
Thus, integrating exponentials over simplexes or, more generally, a
polyhedron formed by the disjoint union of simplexes, will generate
divided differences of the exponential.

We shall also
need two simple preliminary results. Let us use $\PP_n$ to denote
the vector space of polynomials of degree $n$.

\begin{lem} We have
\be
\exp(\mu) \exp[\lambda_0, \ldots, \lambda_m]
= \exp[\lambda_0 + \mu, \ldots, \lambda_m + \mu],
\label{a7d}
\ee
where $\lambda_0, \ldots, \lambda_m$ and $\mu$ can be any complex
numbers. 
\end{lem}

\begin{proof}
Immediate.
\end{proof}

\begin{lem}
Let $f \colon \CC \to \CC$ and let $a_1, \ldots, a_n$ be distinct
nonzero complex numbers. Then
\be
f[0, a_1^2, \ldots, a_n^2]
= g[-a_n, \ldots, -a_1, 0, a_1, \ldots, a_n],
\label{dda0}
\ee
where $g(z) = f(z^2)$, for $z \in \CC$.
\end{lem}

\begin{proof}
Let $p \in \PP_n$ interpolate $f$ at $0, a_1^2, \ldots, a_n^2$.
Then $q(z) := p(z^2)$ is a polynomial of degree $2n$ satisfying
$q(\pm a_j) = p(a_j^2) = f(a_j^2) = g(\pm a_j)$, for $j=0, \ldots, n$,
setting $a_0 = 0$, for convenience. The result then follows from uniqueness
of the interpolating polynomial.

\end{proof}

It is well-known that a divided difference at equally spaced points
can be expressed in a particularly simple form using the forward difference
operator 
\[
\Delta_h f(x) := f(x+h) - f(x),
\]
which we shall need when
demonstrating the equivalence between our moment calculations
and those of \citeasnoun{oshanin} and \citeasnoun{yor}. 
The next proposition is well-known, but we again include
its short proof for the reader's convenience.

\begin{prop} Let $f \colon \RR \to \RR$, let $h$ be any positive
constant and let $n$ be a non-negative integer. Then
\be
f[x, x+h, x+2h, \ldots, x + nh]
= \frac{\Delta_h^n f(x)}{n! h^n}.
\label{dda1}
\ee
\label{dprop5}
\end{prop}

\begin{proof}
It is easily checked that $f[x, x+h] = \Delta_h f(x) / h$. Further,
if we assume \eqref{dda1} for $n-1$, then
the divided difference recurrence relation implies that
\begin{align*}
&f[x, x+h, \ldots, x+nh]\\
&= \frac{f[x+h, \ldots, x+nh] - f[x, x+h, \ldots, x+(n-1)h]}{nh}\\
&= \frac{\Delta_h f[x, \ldots, x+(n-1)h]}{nh}\\
&= \frac{1}{nh} \Delta_h \left( \frac{\Delta_h^{n-1} f(x)}{(n-1)! h^{n-1}}\right)\\
&= \frac{\Delta_h^n f(x)}{n! h^n}.
\end{align*}
Thus the result follows by induction.
\end{proof}

\begin{cor}Let $f \colon \RR \to \RR$ and let $h$ be any positive
constant. Then
\be
f[x, x+h, x+2h, \ldots, x + nh]
= \frac{1}{n! h^n}\sum_{k=0}^n {n \choose k} (-1)^{n-k} f(x + kh).
\label{dda2}
\ee
\label{dcor6}
\end{cor}

\begin{proof}
We define the {\em forward shift} operator
\[
E_h f(x) := f(x+h), \qquad x \in \RR,
\]
and observe
that, by the binomial theorem,
\[
\Delta_h^n f(x) = \left(E_h - 1\right)^n f(x)
= \sum_{k=0}^n {n \choose k} (-1)^{n-k} E_h^k f(x)
= \sum_{k=0}^n {n \choose k} (-1)^{n-k} f(x + kh).
\]
\end{proof}

\begin{cor}
Let $f \colon \RR \to \RR$ and let $h$ be any positive number.
Then
\be
f[-nh, -(n-1)h, \ldots, -h, 0, h, \ldots, nh]
=
\frac{1}{(2n)! h^{2n}} \sum_{k=0}^{2n} {2n \choose k} (-1)^k f((k-n)h).
\label{a7e}
\ee
\label{dcor7}
\end{cor}

\begin{proof}
This is an immediate consequence of Corollary \ref{dcor6}. 
\end{proof}

We shall also need the {\em Leibniz} relation for divided differences
of a product when deriving the recurrence differential equation
for moments.

\begin{thm}[Leibniz] Let $D$ be any subset of $\CC$ containing the distinct points
$z_0, z_1, \ldots, z_n$ and let $v$ and $w$ be complex-valued functions on $D$.
If $u = v \cdot w$, then
\be
u[z_0, \ldots, z_n] = \sum_{k=0}^n v[z_0, \ldots, z_k] w[z_k, \ldots, z_n].
\label{r7}
\ee
\label{dcor9}
\end{thm}

\begin{proof}
Let $p \in \PP_n$ be the unique polynomial interpolant for $u$ 
written in standard Newton form, that is,
\be
p(z) = v[z_0] + v[z_0, z_1](z - z_0) +
\cdots v[z_0, z_1, \ldots, z_n] (z - z_0) \cdots (z - z_{n-1}).
\label{r8}
\ee
We shall let $q \in \PP_n$ be the unique polynomial interpolating $w$, but with
the points chosen in the order $z_n, z_{n-1}, \ldots, z_0$, that is,
\be
q(z) = w[z_n] + w[z_n, z_{n-1}](z - z_n) + \cdots
+ w[z_n, \ldots, z_0] (z - z_n) \cdots (z - z_1).
\label{r9}
\ee
Now their product $p \cdot q$ is a polynomial of degree $2n$. Dividing this
polynomial by $(z - z_0) \cdots (z - z_n)$, we obtain
\[
p(z) q(z) = r(z) + s(z) (z-z_0) \cdots (z-z_n),
\]
where $r \in \PP_n$. We see that $u(z_j) = v(z_j) w(z_j) = p(z_j) q(z_j) = r(z_j)$, 
for $0 \le j \le n$. Hence, by uniqueness of the polynomial interpolant
for $u$ in $\PP_n$, we obtain
\be
r(z) = u[z_0] + \cdots + u[z_0, \ldots, z_n] (z - z_0) \cdots (z - z_n).
\label{r10}
\ee
We obtain \eqref{r7} by equating the coefficients of $z^n$ 
in \eqref{r10} and the product of the
expressions in \eqref{r8} and \eqref{r9}, modulo $(z-z_0) \cdots (z-z_n)$.
\end{proof}


\section{A recurrence relation}

The Feynman--Kac formula (see, for example, \citeasnoun{karatzas})
suggests that the moments $E_n(t) := \EE \left( A(t)^n \right)$ of the
time average should satisfy a certain differential equation, 
which we shall also obtain as an illustration of the divided
difference approach.

\begin{thm}
Let $\{c_n\}_{n=1}^\infty$ be any strictly increasing sequence of
positive numbers and define $e_n : (0,\infty) \to \RR$ by the divided
difference
\be
e_n(t) = \exp[0, c_1 t , \ldots, c_n t], \qquad t > 0, \quad n \ge 0.
\label{r1}
\ee
Then
\be
t e_n^\prime(t) = e_n(t) \left( c_n t - n\right) + e_{n-1}(t),\qquad\hbox{ for } n \ge 1.
\label{r2}
\ee
\label{thm51}
\end{thm}

\begin{proof}
Applying the Hermite--Genocchi formula, we obtain
\be
e_n(t)
=
\int_{K(I_n)} \exp\left(t c^T y\right)\,dy,
\label{r3}
\ee
where $b = \left(c_1, \ldots, c_n\right)^T$, and differentiating \eqref{r3} yields
\be
e_n^\prime(t) = \int_{K(I_n)} \exp(t c^T y) (c^T y)\,dy.
\ee
Now writing $g(s) = s$ and applying Leibniz's formula for divided
differences, we find
\begin{align}
\left(g \cdot \exp\right)[0, c_1 t , \ldots, c_n t]
&= g[0] \exp[0, c_1 t, \ldots, c_n t] + g[0, c_1 t] \exp[c_1 t,
\ldots, c_n t]\nonumber\\
&= \exp[c_1 t, \ldots, c_n t].
\label{r4}
\end{align}
Further, the relation $(g\cdot\exp)^{(n)} = g \cdot\exp + n\exp$ and \eqref{r4} imply
\begin{align}
t e_n^\prime(t)
&= \int_{K(I_n)} \left(g\cdot \exp\right)^{(n)}(t c^T y)\,dy
- n \int_{K(I_n)} \exp(tc^T y)\,dy\nonumber\\
&= \left(g \cdot \exp\right)[0, c_1 t, \ldots, c_n t] - 
n \exp[0, c_1 t, \ldots, c_n t]\nonumber\\
&= \exp[c_1 t, \ldots, c_n t] - n e_n(t).
\label{r5}
\end{align}
However, 
\be
e_n(t) = \frac{\exp[c_1 t, \ldots, c_n t] - e_{n-1}(t)}{c_n t},
\ee
by the divided difference recurrence relation, so that
\be
\exp[c_1 t, \ldots, c_n t]
= c_n t e_n(t) + e_{n-1}(t).
\label{r6}
\ee
Substituting \eqref{r6} in \eqref{r5} provides \eqref{r2}.
\end{proof}

\noindent
The corresponding differential equation for $E_n$ is now immediate.

\begin{cor} The moments satisfy
\be
t E_n^\prime(t) = E_n(t) \left(b_n t - n\right) + E_{n-1}(t), \qquad
\hbox{ for } n \ge 1,
\ee
where $b_n$ is given by \eqref{b1}.
\end{cor}

\begin{proof}
We apply Theorem \ref{thm51} and Theorem \ref{thm_b1}.
\end{proof}

\bibliographystyle{actaagsm}
\bibliography{asian}

\end{document}